\documentclass{elsarticle}
\usepackage{graphicx,amsmath,amsthm,amsfonts,latexsym,amssymb,hyperref,comment,subeqnarray}
 \usepackage{float}
 \usepackage{caption}
\usepackage{hyperref,mathrsfs}
\usepackage[mathscr]{euscript}
\usepackage{colortbl}
\usepackage{booktabs}
\usepackage{todonotes}
\usepackage{algorithmic}
\let\chapter\section               
\usepackage{algorithm}  

\journal{}

\newcommand{\aver}[1]{\left\{\!\!\left\{#1\right\}\!\!\right\}}

\newcommand{\jump}[1]{\left[\!\left[#1\right]\!\right]}
\newcommand{\bigjump}[1]{\left[\!\!\left[#1\right]\!\!\right]}

\newtheorem{remark}{Remark}[section]
\newtheorem{theorem}{Theorem}[section]

\numberwithin{equation}{section}

\begin{document}

\begin{frontmatter}

\title{A Bilinear Partially Penalized Immersed Finite Element Method for Elliptic Interface Problems with Multi-Domains and Triple-Junction Points}

\author{
Yuan Chen\footnote{Department of Statistics, The George Washington University, Washington, DC 20052, USA  (yuan\_chen@gwu.edu). },~
Songming Hou\footnote{Department of Mathematics and Statistics, Louisiana Tech University, Ruston, LA 71272, USA (shou@latech.edu). This aurthor is partially supported by the Walter Koss Professorship fund made available through Louisiana Board of Regents.},~
Xu Zhang\footnote{Department of Mathematics, Oklahoma State University, Stillwater, OK 74078, USA (xzhang@okstate.edu). This author is partially supported by National Science Foundation Grant DMS-1720425}
}

\begin{abstract}
In this article, we introduce a new partially penalized immersed finite element method (IFEM) for solving elliptic interface problems with multi-domains and triple-junction points. We construct new IFE functions on elements intersected with multiple interfaces or with triple-junction points to accommodate interface jump conditions. For non-homogeneous flux jump, we enrich the local approximating spaces by adding up to three local flux basis functions. Numerical experiments are carried out to show that both the Lagrange interpolations and the partial penalized IFEM solutions converge optimally in $L^2$ and $H^1$ norms. 
\end{abstract}

\end{frontmatter}


\section{Introduction}
Elliptic interface problems have received wide attention in the past decades. Numerous numerical methods have been developed to generate accurate and effective approximations for interface problems. Conventional numerical methods, such as the finite element method \cite{chen1998finite}, require the mesh to be aligned with the interface. Another class of numerical methods uses unfitted meshes for solving interface problems. In the finite difference framework, since the pioneering working of immersed boundary method \cite{peskin1977numerical}, the immersed interface method \cite{leveque1994immersed}, matched interface and boundary method \cite{yu2007matched}, cut-cell method \cite{ingram2003developments} have been developed. In the  finite element framework, there have been penalty finite element method \cite{babuvska1970finite}, generalized finite element method \cite{babuvska2012stable}, extended finite element method \cite{fries2010extended}, cut finite element method \cite{2015BurmanClausHansboLarsonMassing}, and  immersed finite element method \cite{li1998immersed}, to name only a few.

The immersed finite element method (IFEM) \cite{li1998immersed,li2003new, li2004immersed,he2008approximation,2017CaoZhangZhang,2019LinSheenZhang} is a class of finite element methods that modify the approximation functions, instead of solution meshes, locally around the interface in order to resolve the interface with unfitted mesh. An IFEM adopts standard finite element basis functions on elements away from interfaces but constructs new piecewise-polynomial basis functions on elements cut through by interfaces. In \cite{PP}, a partially penalized immersed finite element method (PPIFEM) was developed by adding consistency terms and penalty terms to the classical Galerkin formulation on the interface edges. The method is proved to converge optimally in the energy norm  \cite{PP} and $L^2$ norm \cite{guo2018improved}. Moreover,  the PPIFEM was extended to handle non-homogenous flux jump conditions in \cite{ji2018partially}. Gradient recovery for PPIFEM was reported in \cite{guo2017superconvergence}.

Most numerical methods aim at interface problems with two sub-domains. There are a few literatures concerning multi-domain interface problems with possibly a triple-junction point. In \cite{xia2011matched}, the matched interface and boundary method has been developed on multi-domain. In \cite{hou2012numerical}, the authors employed a Petrov-Galerkin type method \cite{hou2010numerical} to solve multi-domain interface problems. Related research papers include \cite{wang2018numerical, wang2015numerical}. As for IFEM for multi-domain interface problems with triple junction, we developed 
a piecewise linear IFEM on triangular meshes \cite{Chen2019}. One of the challenges is the construction of IFEM basis functions due to the complicated geometrical configuration of a multi-domain. We showed that the stiffness matrix is symmetric positive definite and the numerical solutions are reasonably accurate. As indicated in \cite{PP} for two-domain interface problems, we also observed that the convergence rates in $L^2$ and $H^1$ norms deteriorate as the mesh size becomes small.  

In this paper, we develop a new IFEM for solving multi-domain interface problems with non-homogeneous flux jump. We use interface-unfitted rectangular Cartesian meshes, and construct the piecewise bilinear shape functions on multi-interface elements. The PPIFEM is used for enhanced accuracy and stability. Non-homogeneous flux jump conditions are handled by enriching the local approximation spaces, inspired by \cite{He2011Immersed} and \cite{ji2018partially}. It is worthwhile to note that in our new construction on multi-interface elements, the number of restrictions (from nodal-value conditions and interface jump conditions) and degrees of freedom match for all types of interface elements. Therefore, the least-squares approximation is unnecessary for obtaining IFE basis functions as the piecewise-linear IFE spaces in \cite{Chen2019}. Our numerical results indicate that the bilinear PPIFEM for multi-domain interface problems converge optimally in $H^1$, $L^2$ and $L^\infty$ norms without deterioration even for fine meshes. Besides, the new method outperforms the classical linear IFEM \cite{Chen2019} in the accuracy around interfaces. 

The rest of this article is organized as follows. In Section 2, we introduce the multi-domain interface problems and recall some preliminary results. In Section 3, bilinear IFE basis functions and flux jump basis functions are constructed on different types of interface elements. Local and global IFE spaces are defined accordingly.  In Section 4, we present the PPIFEM for multi-domain interface problems with the non-homogeneous flux jump. In Section 5, numerical examples are provided to show the features of our new method. 
A brief conclusion is presented in Section 6.

\section{Interface Problems}
\begin{figure}[htbp]
\centering
  \includegraphics[width=.4\textwidth]{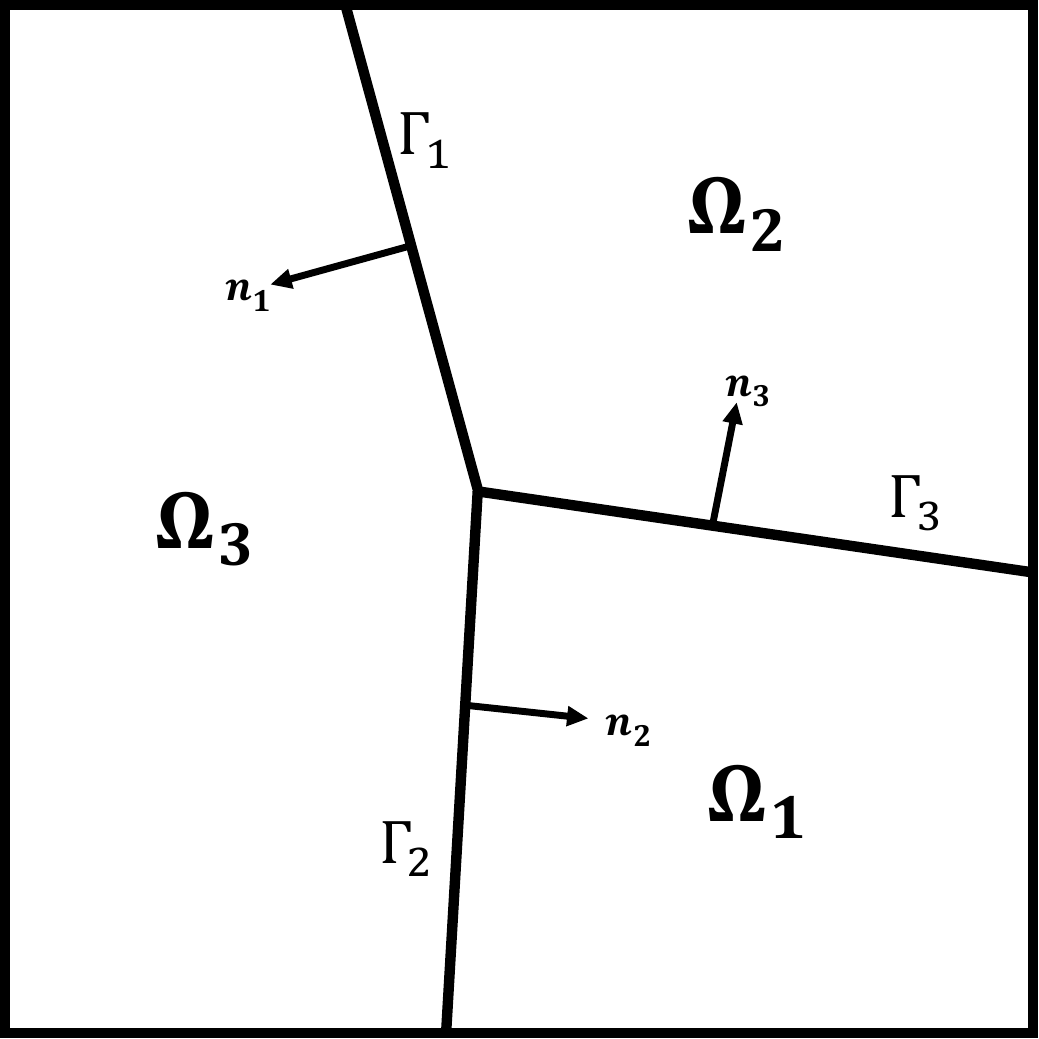}
  \caption{A typical triple-domain domain}
  \label{fig:domain}
\end{figure}

Let $\Omega \subset \mathbb{R}^2$ be an open bounded domain. Assume that $\Omega$ is subdivided by several interfaces  into multiple domains. Without loss of generality, we assume $\Omega$ is divided by three interfaces $\Gamma_1$, $\Gamma_2$, $\Gamma_3$ into three sub-domains $\Omega_1$, $\Omega_2$, $\Omega_3$. We also assume that the boundaries $\partial\Omega_1$, $\partial\Omega_2$, $\partial\Omega_3$ are Lipschitz continuous. See Figure \ref{fig:domain} for an illustration of the geometric setting of the domain and interfaces. 

Consider the following elliptic interface problem:
\begin{equation}\label{pde}
	\begin{split}
		-\nabla\cdot(\beta\nabla u)=f, ~~&\text{in}~ \Omega_1\cup\Omega_2\cup\Omega_3,\\
	u=g,~~&\text{on}~ \partial\Omega.
	\end{split}
\end{equation}
The coefficient function $\beta(\mathbf{x})$ is assumed to be discontinuous across each interface  $\Gamma_i$ where $\mathbf{x} = (x,y)$. For simplicity, we assume that $\beta(\mathbf{x})$ is a piecewise constant function such that 
\begin{equation}\label{jump}
\beta(\mathbf{x}) = \beta_i,~~~~~\text{if}~ \mathbf{x}\in \Omega_i,
\end{equation}
where $\beta_i > 0$.
Across each interface $\Gamma_i$, the solution $u$ is assumed to be continuous, i.e.,
\begin{equation} \label{eq: sol jump}
\jump{u}_{\Gamma_{i}}=0,~~~~\forall i = 1,2,3.
\end{equation}
The normal flux jumps are prescribed as follows:
\begin{equation}\label{eq: flux}
\begin{cases}
\jump{\beta \nabla u \cdot \mathbf{n}}_{\Gamma_{1}}\triangleq(\beta_{3}\nabla u_{3}-\beta_{2}\nabla u_{2})\cdot \mathbf{n}_{1}=b_{1}(x,y)~~~\text{on}~\Gamma_{1},\\
\jump{\beta \nabla u \cdot \mathbf{n}}_{\Gamma_{2}}\triangleq(\beta_{1}\nabla u_{1}-\beta_{3}\nabla u_{3})\cdot \mathbf{n}_{2}=b_{2}(x,y)~~~\text{on}~\Gamma_{2},\\
\jump{\beta \nabla u \cdot \mathbf{n}}_{\Gamma_{3}}\triangleq(\beta_{2}\nabla u_{2}-\beta_{1}\nabla u_{1})\cdot \mathbf{n}_{3}=b_{3}(x,y)~~~\text{on}~\Gamma_{3}.
\end{cases}
\end{equation}
Here, $u_i = u|_{\Omega_i}$, and  $b_i$, $i=1,2,3$, are given functions defined on $\Gamma_i$.  Unit normal vectors of $\Gamma_1$, $\Gamma_2$, $\Gamma_3$ are denoted by $\textbf{n}_1$, $\textbf{n}_2$, $\textbf{n}_3$, respectively. For simplicity, we also use  $g_i = g|_{\partial\Omega_i\cap \partial\Omega}$ from now on.

We use the standard notations for the Sobolev spaces. Define the broken Sobolev space
\[
\tilde H^2(\Omega) = \left\{v \in H^1(\Omega): v|_{\Omega_i}\in H^2(\Omega_i), ~
\jump{\beta \nabla u \cdot \mathbf{n}}_{\Gamma_{i}} = b_i,~i=1,2,3\right\},
\]
equipped with the semi-norm and norm:
\[
|u|_{\tilde H^2(\Omega)}  = \left(\sum_{i=1}^3 |u|_{H^2(\Omega_i)}^2\right)^{1/2},~~~
\|u\|_{\tilde H^2(\Omega)}  = \left(\|u\|_{H^1(\Omega)}^2 +|u|_{\tilde H^2(\Omega)} ^2\right)^{1/2}.
\]
The variational form of \eqref{pde}-\eqref{eq: flux} is to find $u\in H^1(\Omega)$ such that $u=g$ on $\partial \Omega$ and 
\begin{equation}\label{eq: weak form}
a(u,v) :=  (\beta \nabla u,\nabla v) = (f,v) - \sum_{i=1}^3( b_i,v)_{\Gamma_i},~~~\forall v\in H_0^1(\Omega)
\end{equation}
where $(\cdot,\cdot)_\omega$ is the $L^2$ inner product on $\omega\subseteq\Omega$. The subscript $\omega$ is omitted when $\omega = \Omega$. 

For elliptic interface problems, the following regularity result holds \cite{1996BrambleKing}:
\begin{theorem}
Assume that $f\in L^2(\Omega)$, $g\in H^{3/2}(\partial\Omega_i\cap\partial\Omega)$, for $i=1,2,3$. Then the problem \eqref{pde} - \eqref{eq: flux} has a unique solution 
$u\in H^1(\Omega)$ such that for some constant $C>0$,
\begin{equation}
\|u\|_{\tilde H^2(\Omega)} \leq C\left( \|f\|_{L^2(\Omega)}+ \sum_{i=1}^3 \|b_i\|_{H^{1/2}(\Gamma_i)} +  \sum_{i=1}^3 \|g_i\|_{H^{3/2}(\partial\Omega_i\cap\partial\Omega)}\right).
\end{equation}
\end{theorem}

\section{Construction of Triple-Junction Bilinear IFE Spaces}
In this section, we construct the bilinear IFE space on interface elements with multiple interfaces. The construction of classical bilinear IFE functions with one interface is reported in  \cite{Lin2001,He2008}.
\subsection{Element Classification}

\begin{figure}[htbp]
  \includegraphics[width=.49\textwidth]{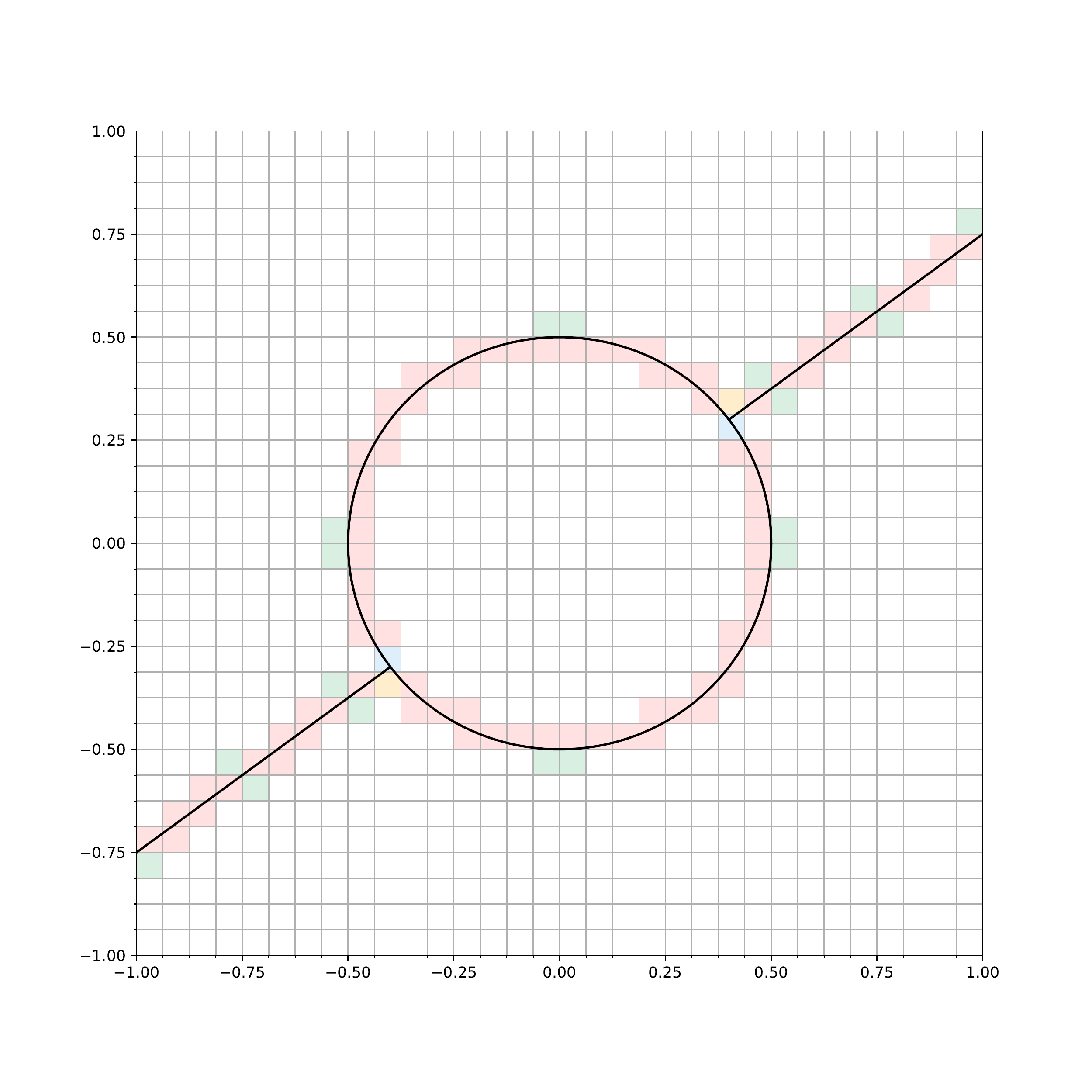}
  \includegraphics[width=.49\textwidth]{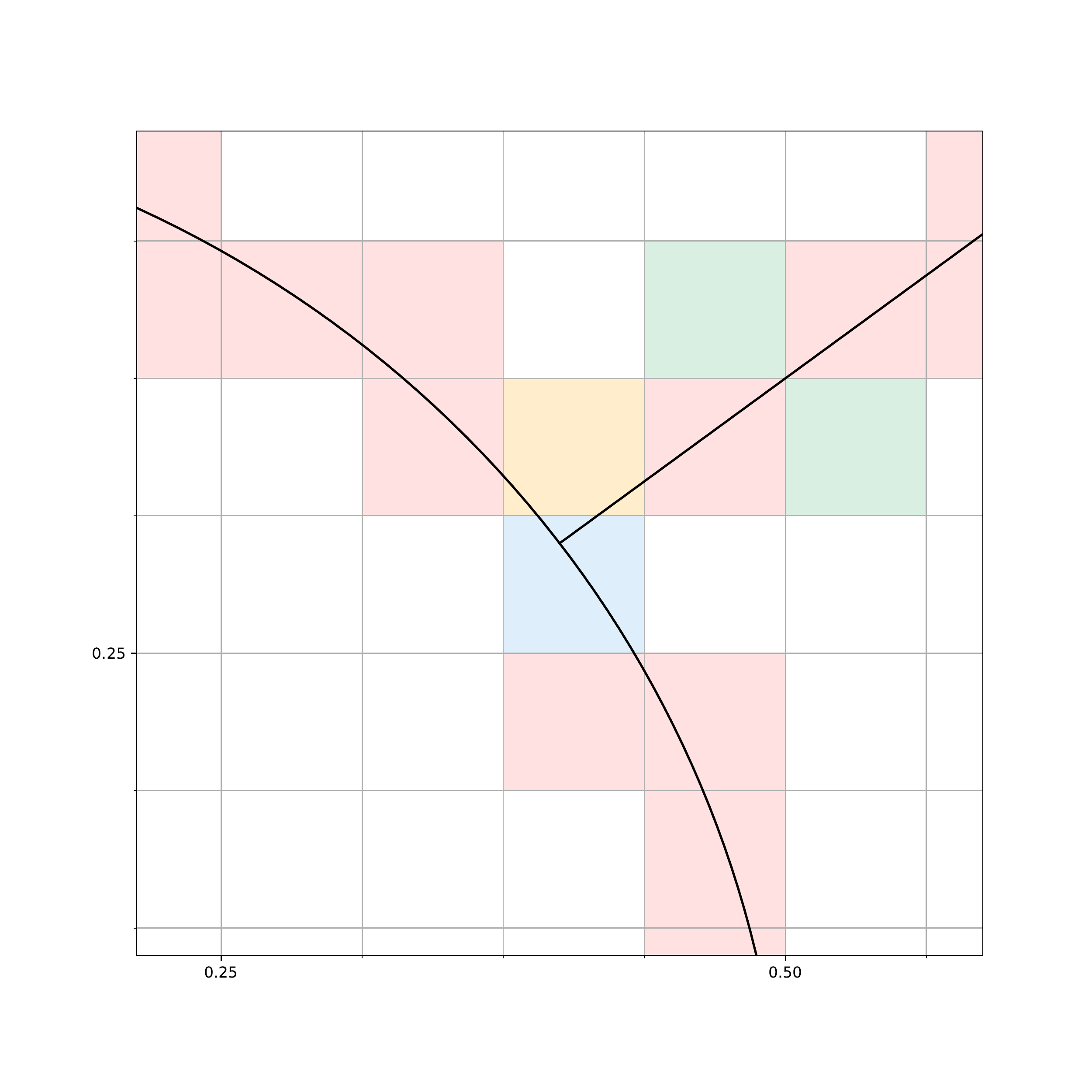}
  \caption{Left: a typical triple-domain domain and Cartesian mesh. Right: a zoom-in of the region around triple point.}
  \label{fig:A typical triple-domain domain}
\end{figure}

\begin{figure}[h]
  \centering
  \includegraphics[width=.7\textwidth,height=.35\textwidth]{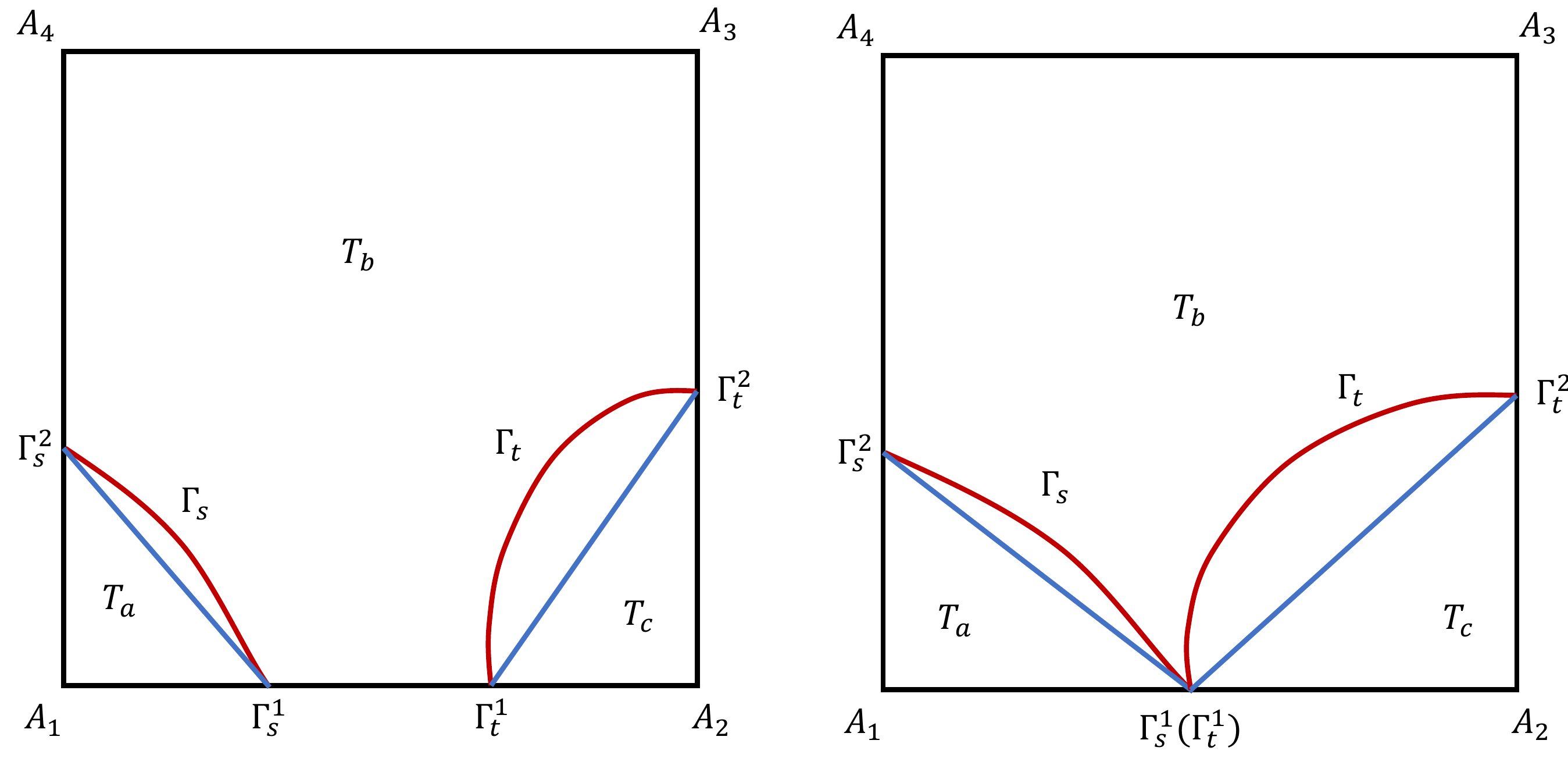}
  \caption{Two typical elements in $\mathcal{T}_{h,2}^i$}
  \label{fig:Two typical elements in 2}
\end{figure}

\begin{figure}[h]
  \centering
  \includegraphics[width=.7\textwidth,height=.35\textwidth]{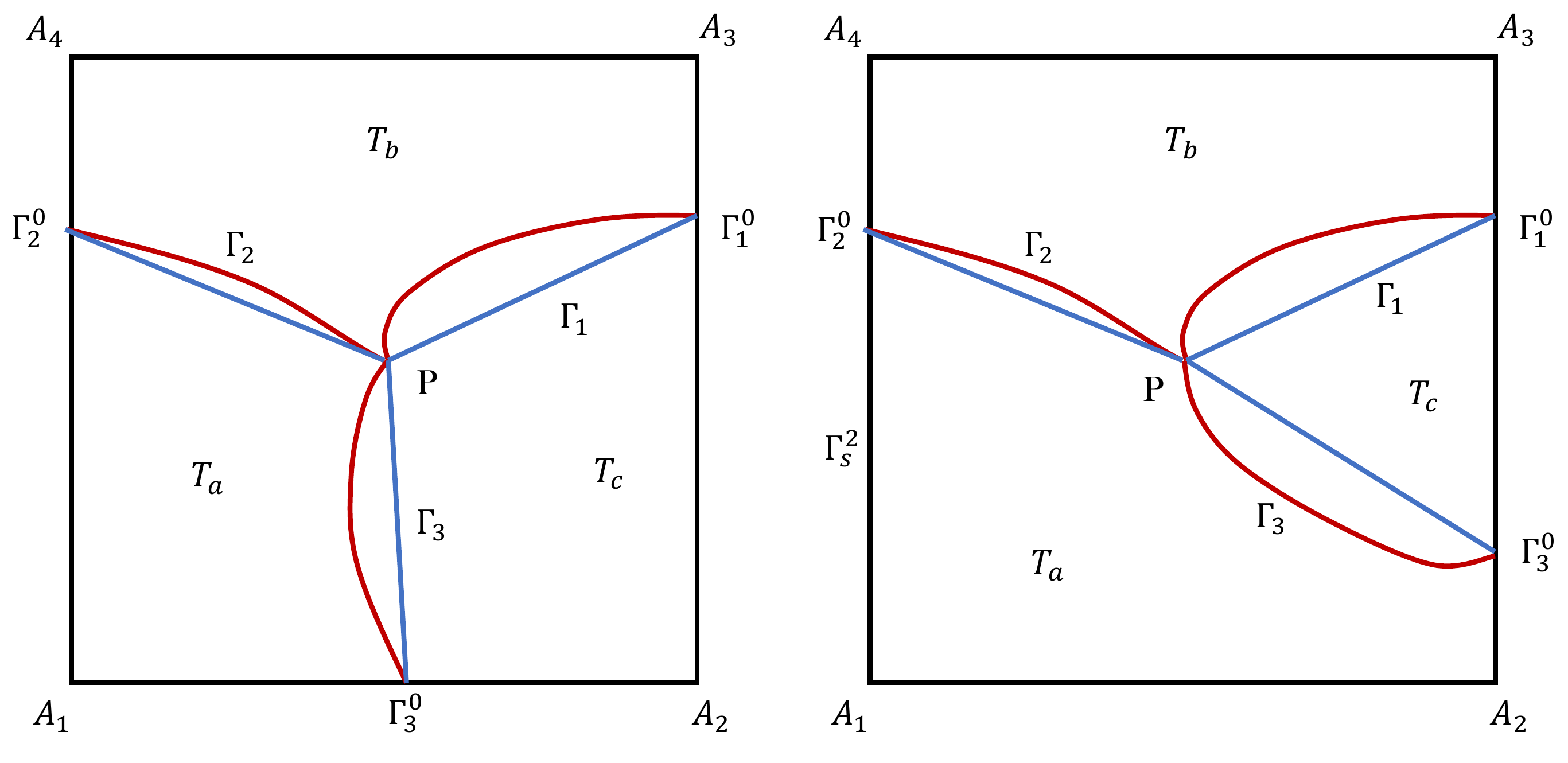}
  \caption{Two typical elements in $\mathcal{T}_{h,3}^i$}
  \label{fig:Two typical elements in 3}
\end{figure}

Let $\mathcal{T}_h$ be an interface-independent Cartesian mesh of $\Omega$ where $h$ denotes the mesh size. 
We classify all rectangular elements according to the number of interfaces intersecting the interior of an element. If the interior of an element $T\in \mathcal{T}_h$ is not intersected with any interface $\Gamma_i$, then the element $T$ is called a non-interface (or a regular) element. The collection of all regular elements are denoted by $\mathcal{T}_{h}^n$. 
If the interior of an element $T\in \mathcal{T}_h$ is cut through by at least one interface $\Gamma_i$, then the element $T$ is called an interface element. The collection of these elements is denoted by $\mathcal{T}_{h}^i$. For multi-interface problems, we assume that the interface elements satisfy the following hypotheses:
\begin{description}
\item[(H1)] Multiple interfaces can only intersect with an element at no more than three edges. 
\item[(H2)] Each edge of an interface element can only intersect with interfaces at no more than two points unless the edge is part of the interface. 
\end{description}
These hypotheses rule out the case that interfaces intersect an element at three points but on the same edge. Furthermore, we distinguish the collection of interface elements $\mathcal{T}_{h}^i$ by three sub-categories based on the number of interfaces inside an element. For a triple-domain interface problem stated above, there are at most three interfaces inside an element. 
The set $\mathcal{T}_{h,1}^i$ contains elements intersected with only one interface. The set $\mathcal{T}_{h,2}^i$ includes all elements intersecting two interfaces (see Figure \ref{fig:Two typical elements in 2} for an illustration). The set $\mathcal{T}_{h,3}^i$ contain the elements intersecting with all three interfaces and hence the triple-junction points must be inside these elements (see Figure \ref{fig:Two typical elements in 3} for an illustration).

Now all elements in a Cartesian mesh must be in one of the following classes.
\begin{equation}
\mathcal{T}_{h}=\mathcal{T}_{h}^n \cup \mathcal{T}_{h}^i  = \mathcal{T}_{h}^n \cup \mathcal{T}_{h,1}^i \cup \mathcal{T}_{h,2}^i \cup \mathcal{T}_{h,3}^i.
\end{equation}

%
%
%

Figure \ref{fig:A typical triple-domain domain} illustrates the above classification of elements by color. Elements in white denote regular elements. Elements in green denote an interface touches one of the nodes of the element or overlaps with an edge of the element. Elements in both of these two colors represent elements in $\mathcal{T}_{h}^n$. Elements marked in red, orange, and blue represent elements in $\mathcal{T}_{h,1}^i$, $\mathcal{T}_{h,2}^i$, and $\mathcal{T}_{h,3}^i$, respectively. 
\begin{remark}
If the triple-junction point is on an edge of some element (see the right plot of Figure \ref{fig:Two typical elements in 2}), we count this element in $\mathcal{T}_{h,2}^i$. 
\end{remark}

\subsection{Local IFE spaces on $T \in \mathcal{T}_{h}^n$ and $T \in \mathcal{T}_{h,1}^i$}
On a non-interface element $T\in \mathcal{T}_{h}^n$, we adopt the standard bilinear ($\mathbb{Q}_1$) nodal basis functions. The local finite element space on $T$ is \[S_h(T) = \text{span}\{1,x,y,xy\}.\]
If an element $T \in \mathcal{T}_{h,1}^i$, then $T$ intersects with only one interface at two different points at two edges. We refer the readers to the construction of the classical bilinear IFE basis functions in \cite{Lin2001,He2008}. For the non-homogeneous flux jump, the local space is enlarged by adding one more basis function that vanishes at all nodes, continuous over the element, but with a unit flux jump. For more details, we refer readers to \cite{He2011Immersed}. Our main focus is to construct the local approximating spaces on elements in $\mathcal{T}_{h,2}^i$ and $\mathcal{T}_{h,3}^i$.

\subsection{Local IFE spaces on $T \in \mathcal{T}_{h,2}^i$}
In this case, there are two interfaces inside an element $T \in \mathcal{T}_{h,2}^i$. As shown in Figure \ref{fig:Two typical elements in 2}, the triple-junction point is disjoint of the interior of the element (left plot of Figure \ref{fig:Two typical elements in 2}) or lie on the edge of the element (right plot of Figure \ref{fig:Two typical elements in 2}). Let the curved interfaces be approximated by the line segments  $\Gamma_s,~\Gamma_t$, and they split the element $T$ into three polygons, denoted by $T_a,~T_b,~T_c$. The intersection points are denoted by $\Gamma_s^1$, $\Gamma_s^2$, $\Gamma_t^1$, and $\Gamma_t^2$. See Figure \ref{fig:Two typical elements in 2} for the geometrical configuration.

\begin{figure}[h]
  \centering
  \includegraphics[width=.49\textwidth]{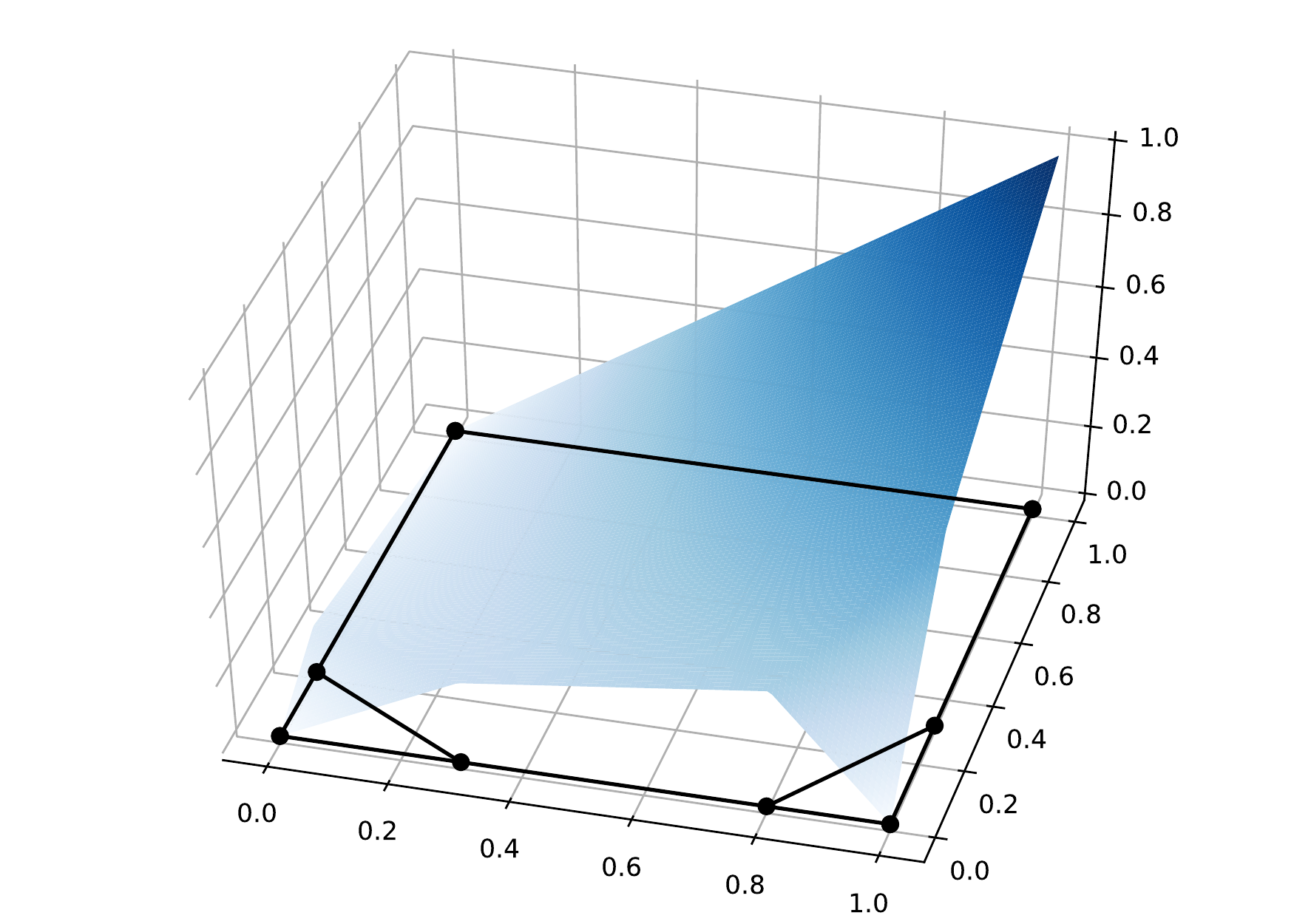}
  \includegraphics[width=.49\textwidth]{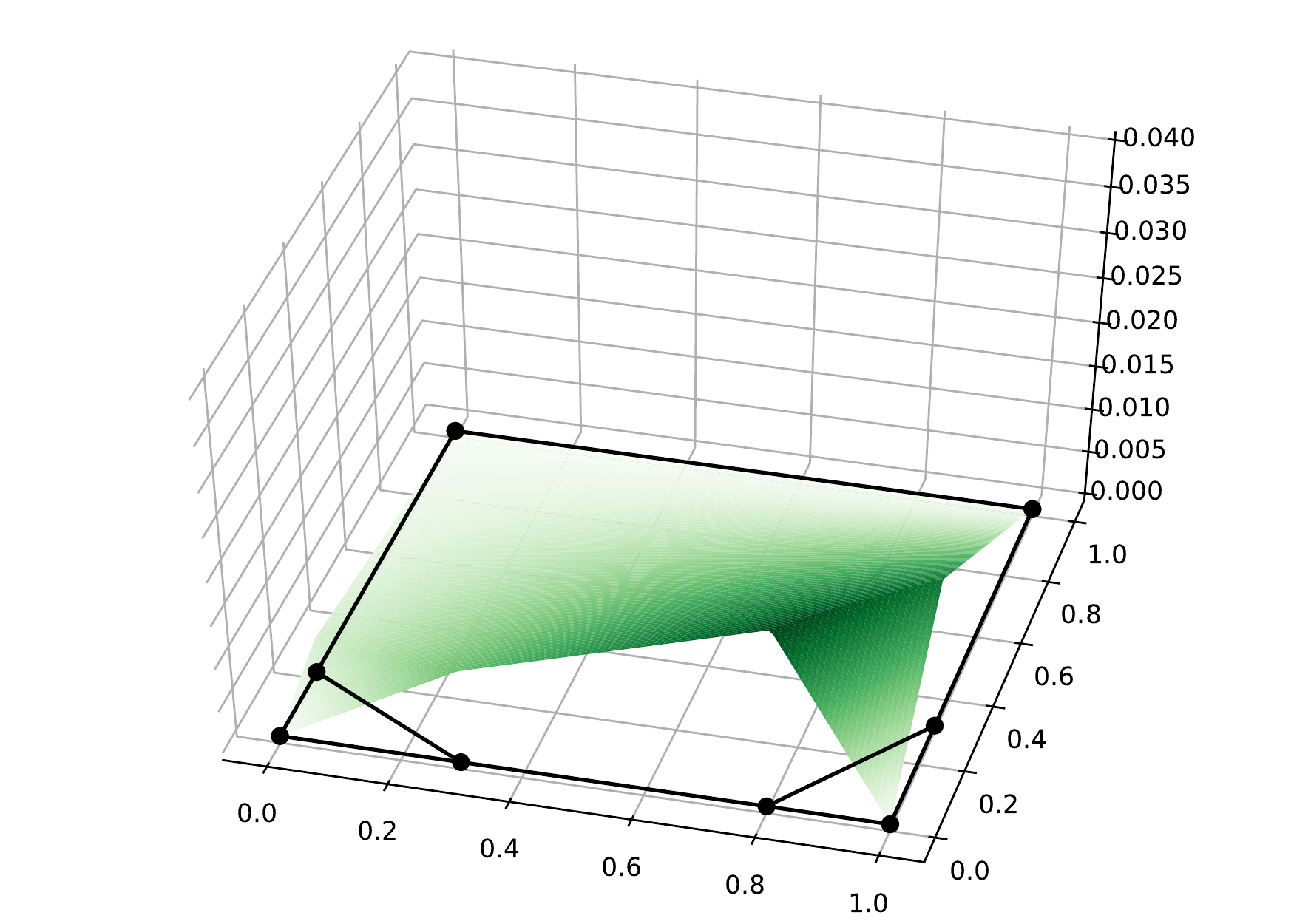}
  \caption{Basis function(left) and flux jump function(right) of elements in $\mathcal{T}_{h,2}^i$}
  \label{fig:Plot in 4}
\end{figure}

We define four piecewise bilinear basis functions on an element as follows
\begin{equation}
\phi_{i,T}(\mathbf{x})=
	\begin{cases}
		\phi_{i,T}^a(\mathbf{x})=a_1+b_1x+c_1y+d_1xy,~~\text{if}~\mathbf{x}\in T_{a}, \\
		\phi_{i,T}^b(\mathbf{x})=a_2+b_2x+c_2y+d_2xy,~~\text{if}~\mathbf{x}\in T_{b}, \\
		\phi_{i,T}^c(\mathbf{x})=a_3+b_3x+c_3y+d_3xy,~~\text{if}~\mathbf{x}\in T_{c}.
	\end{cases}~~i=1,2,3,4.
\label{equ:basis equation points4}
\end{equation}
The following nodal-value conditions and the interface jump conditions are used to determine the coefficients $a_j$, $b_j$, $c_j$, and $d_j$, $(j=1,2,3)$:
\begin{itemize}
\item four nodal-value conditions
\begin{equation}
\phi_{i,T}(A_j)=\delta_{ij},~~~\forall i,j = 1,2,3,4.
\label{eq: nodal value 2-1}
\end{equation}
\item six continuity conditions of the basis functions
\begin{equation}\label{eq: sol continuity 2-1}
\begin{split}
\phi_{i,T}^a(\Gamma_s^1)=\phi_{i,T}^b(\Gamma_s^1),~~
\phi_{i,T}^a(\Gamma_s^2)=\phi_{i,T}^b(\Gamma_s^2),~~
d_1 = d_2,\\
\phi_{i,T}^b(\Gamma_t^1)=\phi_{i,T}^c(\Gamma_t^1),~~
\phi_{i,T}^b(\Gamma_t^2)=\phi_{i,T}^c(\Gamma_t^2),~~
d_2=d_3.\end{split}
\end{equation}
\item two conditions of normal flux continuity
\begin{equation}
\int_{\Gamma_s\cap T}\jump{\beta \nabla \phi_{i,T}\cdot\mathbf{n}}_{\Gamma_{s}}ds=0,~~~~
\int_{\Gamma_t\cap T}\jump{\beta \nabla \phi_{i,T}\cdot\mathbf{n}}_{\Gamma_{t}}ds=0.
\label{eq: flux continuity 2-1}
\end{equation}
\end{itemize}
The twelve conditions in \eqref{eq: nodal value 2-1} - \eqref{eq: flux continuity 2-1} are used to determine the twelve coefficients $a_i$, $b_i$, $c_i$, and $d_i$, with $i=1,2,3$ in a bilinear IFE local basis function $\phi_{i,T}$,  $i=1,2,3,4$. 
The left plot of Figure \ref{fig:Plot in 4} exemplifies a typical bilinear IFE basis function on $T\in \mathcal{T}_{h,2}^i$. The local bilinear IFE space is defined as 
$S_h(T) = \text{span}\{\phi_{i,T}: i=1,2,3,4\}$.

To accommodate the non-homogeneous jump of flux \eqref{eq: flux}, we enrich the local IFE space by adding two local functions $\phi_{T,J}^k$, $k=1,2$. These flux basis functions are also defined as piecewise bilinear polynomials 
\begin{equation}
\phi_{T,J}^{k}(\mathbf{x})=
	\begin{cases}
		\phi_{T,J}^{k,a}(\mathbf{x})=a_1+b_1x+c_1y+d_1xy,~~\text{if}~\mathbf{x}\in T_{a}, \\
		\phi_{T,J}^{k,b}(\mathbf{x})=a_2+b_2x+c_2y+d_2xy,~~\text{if}~\mathbf{x}\in T_{b}, \\
		\phi_{T,J}^{k,c}(\mathbf{x})=a_3+b_3x+c_3y+d_3xy,~~\text{if}~\mathbf{x}\in T_{c}.
	\end{cases}
	k=s,t.
\label{eq: flux basis 3 pts}
\end{equation}
The flux jump basis functions $\phi_{T,J}^s$ and $\phi_{T,J}^t$ satisfy the following conditions
\begin{itemize}
\item They vanish at four nodes $A_j$ of the element $T$:
\begin{equation}\label{}
\phi_{T,J}^k(A_j)=0, ~~~\forall j = 1,2,3,4,~~k = s,t.
\end{equation}
\item They are continuous across the endpoints of interfaces $\Gamma_i\cap T$:
\begin{equation}\label{eq: continuity flux bas}
\begin{split}
\begin{split}
\phi_{T,J}^{k,a}(\Gamma_s^1)=\phi_{T,J}^{k,b}(\Gamma_s^1),~~
\phi_{T,J}^{k,a}(\Gamma_s^2)=\phi_{T,J}^{k,b}(\Gamma_s^2),~~
d_1 = d_2,\\
\phi_{T,J}^{k,a}(\Gamma_t^1)=\phi_{T,J}^{k,b}(\Gamma_t^1),~~
\phi_{T,J}^{k,a}(\Gamma_t^2)=\phi_{T,J}^{k,b}(\Gamma_t^2),~~
d_2=d_3.\end{split}
\end{split}
\end{equation}
\item They have a unit flux jump on the corresponding interfaces. 
\begin{equation}\label{eq: flux jump flux bas}
\int_{\Gamma_{k'}\cap T}\bigjump{\beta \nabla \phi_{T,J}^k\cdot\mathbf{n}}_{\Gamma_{k'}}ds=\delta_{kk'},~~~\forall k' = s,t. 
\end{equation}
\end{itemize}
A flux basis function of this type is shown in the right plot of Figure \ref{fig:Plot in 4}. The local flux IFE space $S_h^J(T) = \text{span}\{\phi_{T,J}^s, \phi_{T,J}^t\}$
and the enriched local IFE space on  $T$ is $\tilde S_h(T) = S_h(T) \cup S_h^J(T)$. The Lagrange type interpolation $\mathcal{I}_{h,T}: \tilde H^2(T) \to \tilde S_h(T)$ is as follows
\begin{equation}
\mathcal{I}_{h,T}u(\mathbf{x})=\sum_{i=1}^3u(A_i)\phi_{i,T}(\mathbf{x})+ \sum_{k=s,t} q_T^k\phi_{T,J}^k (\mathbf{x}),
\end{equation}
where 
\begin{equation}
	q_T^k=\int_{\Gamma_k\cap T} \jump{\beta \nabla u\cdot \mathbf{n}}_{\Gamma_{k}}ds.
\label{equ:equation c}
\end{equation}


\subsection{Local IFE spaces on $T \in \mathcal{T}_{h,3}^i$}
In this case, an element $T \in \mathcal{T}_{h,3}^i$ intersects all three interfaces at three points. In general, the triple-junction point, denoted by $P$, is inside this element. See Figure \ref{fig:Two typical elements in 3} for an illustration. 

We use the line approximations of the curved interfaces by $\Gamma_1,~\Gamma_2,~\Gamma_3$, which split the whole rectangle into three polygons, denoted by $T_a,~T_b,~T_c$ (see figure \ref{fig:Two typical elements in 3}). The intersection points with each interface $\Gamma_1,~\Gamma_2,~\Gamma_3$ are denoted by $\Gamma_1^0$, $\Gamma_2^0$, and $\Gamma_3^0$, respectively.

\begin{figure}[h]
  \centering
  \includegraphics[width=.49\textwidth]{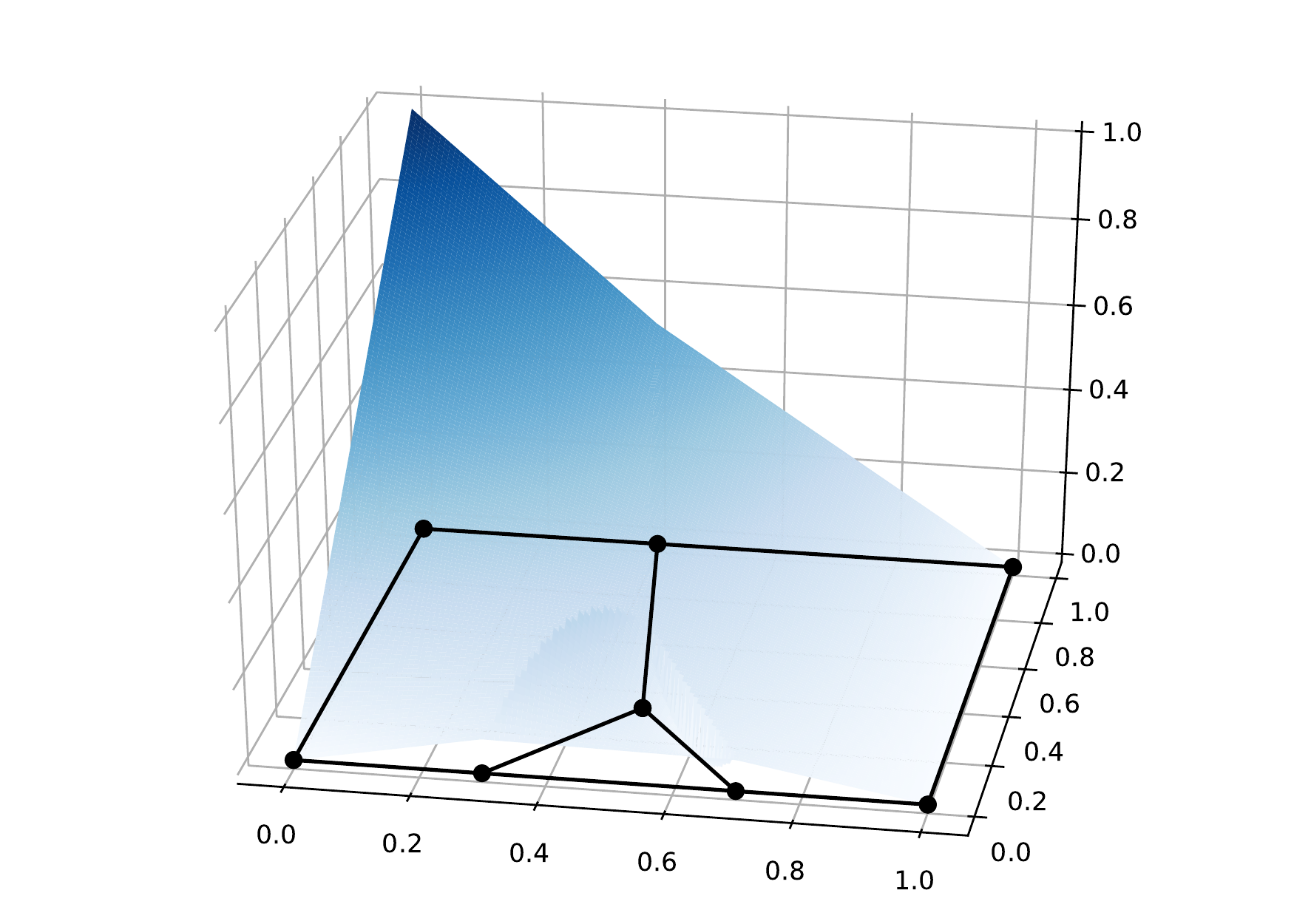}
  \includegraphics[width=.49\textwidth]{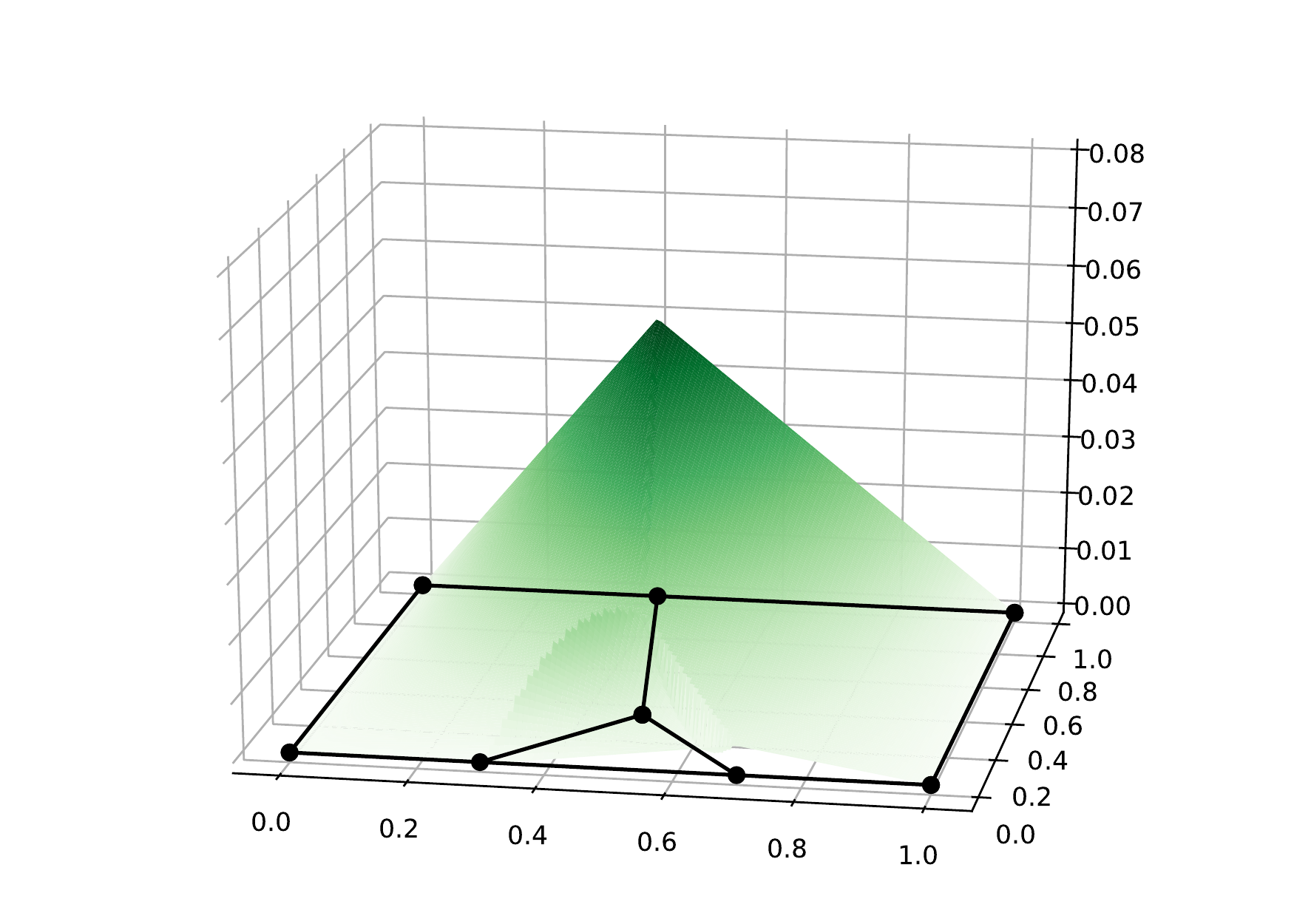}
  \caption{A bilinear IFE local basis function(left) and a flux basis function (right) of elements in $\mathcal{T}_{h,3}^i$}
  \label{fig:Plot in 3}
\end{figure}

Piecewise bilinear IFE basis functions are constructed on the element $T$ as follows
\begin{equation}
\phi_{i,T}(\mathbf{x})=
	\begin{cases}
		\phi_{i,T}^a(\mathbf{x})=a_1+b_1x+c_1y+d_1xy,~~\text{if}~\mathbf{x}\in T_{a}, \\
		\phi_{i,T}^b(\mathbf{x})=a_2+b_2x+c_2y+d_2xy,~~\text{if}~\mathbf{x}\in T_{b}, \\
		\phi_{i,T}^c(\mathbf{x})=a_3+b_3x+c_3y+d_3xy,~~\text{if}~\mathbf{x}\in T_{c}.
	\end{cases}
	i = 1,2,3,4.
\label{equ:basis equation points3}
\end{equation}
As before, the twelve coefficients are determined by 
\begin{itemize}
\item four-nodal value conditions
\begin{equation}
\phi_{i,T}(A_j)=\delta_{ij},~~~\forall i,j = 1,2,3,4.
\label{eq: nodal value 2}
\end{equation}
\item five continuity conditions of the basis functions
\begin{equation}\label{eq: sol continuity 2}
\begin{split}
&\phi_{i,T}^b(\Gamma_1^0)=\phi_{i,T}^c(\Gamma_1^0),~~
\phi_{i,T}^a(\Gamma_2^0)=\phi_{i,T}^b(\Gamma_2^0),~~
\phi_{i,T}^a(\Gamma_3^0)=\phi_{i,T}^c(\Gamma_3^0),\\
&\phi_{i,T}^a(P)=\phi_{i,T}^b(P) ,~~
\phi_{i,T}^b(P) = \phi_{i,T}^c(P). 
\end{split}
\end{equation}
\item three conditions of normal flux continuity
\begin{equation}\label{eq: flux continuity 2}
\int_{\Gamma_k\cap T}[\beta \nabla \phi_{i,T}\cdot\mathbf{n}]_{\Gamma_{k}}ds=0,~~~k=1,2,3.
\end{equation}
\end{itemize}
A typical bilinear IFE basis function is plotted in the left plot of Figure \ref{fig:Plot in 3}. The local bilinear IFE space is defined as 
$S_h(T) = \text{span}\{\phi_{i,T}: i=1,2,3,4\}$.

For non-homogeneous jump of flux \eqref{eq: flux}, we enrich the local IFE space by adding three local functions $\phi_{T,J}^k$, $k=1,2,3$ as follows
\begin{equation}
\phi_{T,J}^{k}(\mathbf{x})=
	\begin{cases}
		\phi_{T,J}^{k,a}(\mathbf{x})=a_1+b_1x+c_1y+d_1xy,~~\text{if}~\mathbf{x}\in T_{a}, \\
		\phi_{T,J}^{k,b}(\mathbf{x})=a_2+b_2x+c_2y+d_2xy,~~\text{if}~\mathbf{x}\in T_{b}, \\
		\phi_{T,J}^{k,c}(\mathbf{x})=a_3+b_3x+c_3y+d_3xy,~~\text{if}~\mathbf{x}\in T_{c}.
	\end{cases}
	k=1,2,3,
\label{eq: flux basis 3 pts}
\end{equation}
such that for $k=1,2,3$ the following conditions hold
\begin{itemize}
\item vanishing at four nodes $A_j$ of the element $T$:
\begin{equation}\label{}
\phi_{T,J}^k(A_j)=0, ~~~\forall j = 1,2,3,4.
\end{equation}
\item continuous across the endpoints of interfaces $\Gamma_i\cap T$:
\begin{equation}\label{eq: continuity flux bas}
\begin{split}
&\phi_{T,J}^{k,b}(\Gamma_1^0)=\phi_{T,J}^{k,c}(\Gamma_1^0),~~
\phi_{T,J}^{k,a}(\Gamma_2^0)=\phi_{T,J}^{k,b}(\Gamma_2^0),~~
\phi_{T,J}^{k,a}(\Gamma_3^0)=\phi_{T,J}^{k,c}(\Gamma_3^0),\\
&\phi_{T,J}^{k,a}(P)=\phi_{T,J}^{k,b}(P) ,~~
\phi_{T,J}^{k,b}(P) = \phi_{T,J}^{k,c}(P). 
\end{split}
\end{equation}
\item unit flux jump condition:
\begin{equation}\label{eq: flux jump flux bas}
\int_{\Gamma_j\cap T}[\beta \nabla \phi_{T,J}^k\cdot\mathbf{n}]_{\Gamma_{j}}ds=\delta_{kj},~~~\forall j = 1,2,3. 
\end{equation}
\end{itemize}
A flux basis function of this type is shown in the right plot of Figure \ref{fig:Plot in 3}.  The local flux space is defined by $S_h^J(T) = \text{span}\{\phi_{T,J}^1, \phi_{T,J}^2, \phi_{T,J}^3\}$. As before, the enriched local IFE space on $T$ is $\tilde S_h(T) = S_h(T) \cup S_h^J(T)$. The Lagrange type interpolation $\mathcal{I}_{h,T}: \tilde H^2(T) \to \tilde S_h(T)$ is as follows
\begin{equation}
\mathcal{I}_{h,T}u(\mathbf{x})=\sum_{i=1}^4u(A_i)\phi_{i,T}(\mathbf{x})+ \sum_{k=1}^3 q_T^k\phi_{T,J}^k (\mathbf{x}),
\end{equation}
where arguments $q_T^k$, $k=1,2,3$, are integrations of flux jump along the interface, i.e.,
\begin{equation}
	q_T^k=\int_{\Gamma_k\cap T} [\beta \nabla u\cdot \mathbf{n}]_{\Gamma_{k}}ds.
\label{equ:equation b}
\end{equation}
\begin{remark}
On an element $T \in \mathcal{T}_{h,2}^i$, the bilinear IFE basis functions degenerate to standard bilinear IFE basis functions in \cite{He2008} when an interface moves out of the element. This can be proved following the same idea in \cite{2019LinSheenZhang}.
\end{remark}
\begin{remark}\label{rem: 3.3}
Compared with the linear IFE basis functions of the triple-junction element in \cite{Chen2019}, the bilinear IFE basis functions and flux basis functions on rectangular element ensure the nodal value conditions precisely. Inside the triple-junction element, the continuity is guaranteed at the intersection points and triple junction point.
\end{remark}

%

\subsection{Global bilinear IFE space}
Now we define the global bilinear  IFE space and flux IFE space. Let $\mathcal{N}_h$ be the collection of all interior nodes of the mesh $\mathcal{T}_h$. 
Associated with each node $\mathbf{x}_i\in \mathcal{N}_h$, we define a global bilinear IFE function $\Phi_i$ such that
\begin{description}
\item[(i)] $\Phi_i(x_j) = \delta_{ij}$, $\forall i,j \in \mathcal{N}_h$.
\item[(ii)] $\Phi_i|_T\in S_h(T)$,~$\forall T\in \mathcal{T}_h$,
\end{description}
where $S_h(T)$ is the local IFE space defined in previous subsections. Correspondingly, the global bilinear IFE space is defined as 
\[
S_h(\Omega)=\text{span}\{\Phi_i: i=1,2,...,|\mathcal{N}_h|\}.
\]
For each local flux basis $\phi_{T,J}^k$, a zero-extension to the whole domain will yield the corresponding global flux basis function $\Phi_{T,J}^k$:
\begin{equation}
    \Phi_{T,J}^k(\mathbf{x})=
	\begin{cases}
		\phi_{T,J}^k(\mathbf{x}),&\text{if}~~\mathbf{x} \in T,\\
		0,                     &\text{if}~~\mathbf{x} \notin T.
	\end{cases}
\end{equation}
Here, the value $k$ depends on the type of interface element. In particular, $k=1$ if $T \in \mathcal{T}_{h,1}^i$, $k=1$, $k=1,2$ if $T \in \mathcal{T}_{h,2}^i$, and $k=1,2,3$ if $T \in \mathcal{T}_{h,3}^i$.
The flux jump IFE space is defined as 
\[S_h^J(\Omega)=\text{span}\{\Phi_{T,J}^k: \forall T\in\mathcal{T}_h^i\},\]
and the enriched global bilinear IFE space is $\tilde S_h(\Omega) = S_h(\Omega)\cup S_h^J(\Omega)$. 

The Lagrange interpolation operator $\mathcal{I}_h:  \tilde H^2(\Omega) \to \tilde S_h$ is
\begin{equation}\label{eq: interpolation}
	\mathcal{I}_h u (\mathbf{x}) = \sum_{j \in \mathcal{N}_h}u(\mathbf{x}_j)\Phi_j(\mathbf{x})+\sum_{T \in \mathcal{T}^{i}_h}\sum_{k=1}^{d_k}q^k_T\Phi_{T,J}^k(\mathbf{x}),
\end{equation}
where $d_k = 1$ for $T\in \mathcal{T}^{i}_{h,2}$, $d_k = 2$ for $T\in \mathcal{T}^{i}_{h,2}$, and $d_k = 3$ for $T\in \mathcal{T}^{i}_{h,3}$. The value of $q^k_T$ is given in \eqref{equ:equation c} and \eqref{equ:equation b}.


%

\section{PPIFEM for Triple-Junction Interface Problems}
Let $\mathcal{E}_h$ be the collection of edges on the mesh $\mathcal{T}_h$. Let $\mathring{\mathcal{E}_h}$ and $\mathcal{E}_h^b$ be the set interior edges and boundary edges, respectively. If an edge $e\in \mathcal{E}_h$ intersects with the interface, we call it an interface edge, otherwise, it is called a non-interface edge. Denote by $\mathcal{E}_h^i$ and $\mathcal{E}_h^n$ the set of interface edges and non-interface edges, respectively. 

For any interior edge $e\in\mathring{\mathcal{E}_h}$, we denote the two elements sharing the edge $e$ by $T_{e,1}$ and $T_{e,2}$. For a function $u$ defined on $T_{e,1} \cup T_{e,2}$, we define the average and jump of $u$ across the edge $e$ respectively by 
\[\aver{u}_e=\frac{1}{2}((u|_{T_{e,1}})|_e+(u|_{T_{e,2}})_e),~~~
\jump{u}_e=(u|_{T_{e,1}})|_e-(u|_{T_{e,2}})_e.\]
For $e\in\mathcal{E}_h^b$, we define $\aver{u}_e=\jump{u}_e=u|_e$.

The bilinear IFE approximation is to find $\tilde u_h\in \tilde S_h(\Omega)$ in the following form
\begin{equation}
	\tilde u_h= u_h + u_h^J \triangleq \sum_{j \in \mathcal{N}_h}u_j\Phi_j+\sum_{T \in \mathcal{T}^{i}_h}\sum_{k=1}^3q^k_T\Phi_{T,J}^k,
\end{equation}
where $u_h\in S_h$ is the unknown function, and $u_h^J\in S_h^J$ can be explicitly constructed. We use the PPIFEM \cite{PP} for solving the nonhomogeneous flux triple-junction interface problem: Find $u_h\in S_h$ such that
\begin{equation}
a_h(u_h,v_h) = (f,v_h)-a_h(u_h^J,v_h)+\sum_{i=1}^3(b_i,v_h)_{\Gamma_i},~~~~\forall v_h\in S_h.
\end{equation}
where the bilinear form
\begin{equation}
	\begin{split}
	    a_h(u_h,v_h)&=\sum_{K\in \mathcal{T}_h}\int_{K}\beta\nabla v_h\cdot\nabla u_h d\mathbf{x}-\sum_{e\in \mathring{\mathcal{E}_h^i}}\int_e
	    \aver{\beta\nabla u_h\cdot\textbf{n}_e}\jump{v_h}ds\\
		&+\epsilon \sum_{e\in \mathcal{E}_h^i}\int_e \aver{\beta\nabla v_h\cdot\textbf{n}_e}\jump{u_h}ds+\sum_{e\in \mathcal{E}_h^i}\int_e\frac{\sigma_e}{|e|}\jump{v_h}\jump{u_h}ds.\\
	\end{split}
	\label{equ:6}
\end{equation}
Here, $\epsilon$ has the following choices. When $\epsilon=-1$, the scheme is called symmetric PPIFEM. When $\epsilon=0, 1$, the scheme is called incomplete {PPIFEM} and non-symmetric PPIFEM, respectively. The stability parameter $\sigma_e>0$.

\section{Numerical Examples}
In this section, we test our new bilinear approximation in two numerical examples. The numerical experiments are carried out on Cartesian meshes with $N\times N$ rectangular elements, starting from $N=16$. We consider the following three norms of the numerical solution:
\[\|u_h-u\|_{L^\infty} = \max_{(x,y)\in \mathcal{N}_h}|u_h(x,y)-u(x,y)|,\]
\[\|u_h-u\|_{L^2} = \int_\Omega |u_h(x,y)-u(x,y)|^2dxdy,\]
\[|u_h-u|_{H^1} = \int_\Omega|\nabla u_h(x,y)-\nabla u(x,y)|^2dxdy,\]
where $\mathcal{N}_h$ denotes the set of all nodes in the mesh $\mathcal{T}_h$.
To avoid redundance, we only report errors of symmetric PPIFEM, since the results from the nonsymmetric and incomplete PPIFEM are similar.

\subsection{Example 1 (straight-line interface)}
\begin{figure}[htbp]
  \includegraphics[width=.49\textwidth]{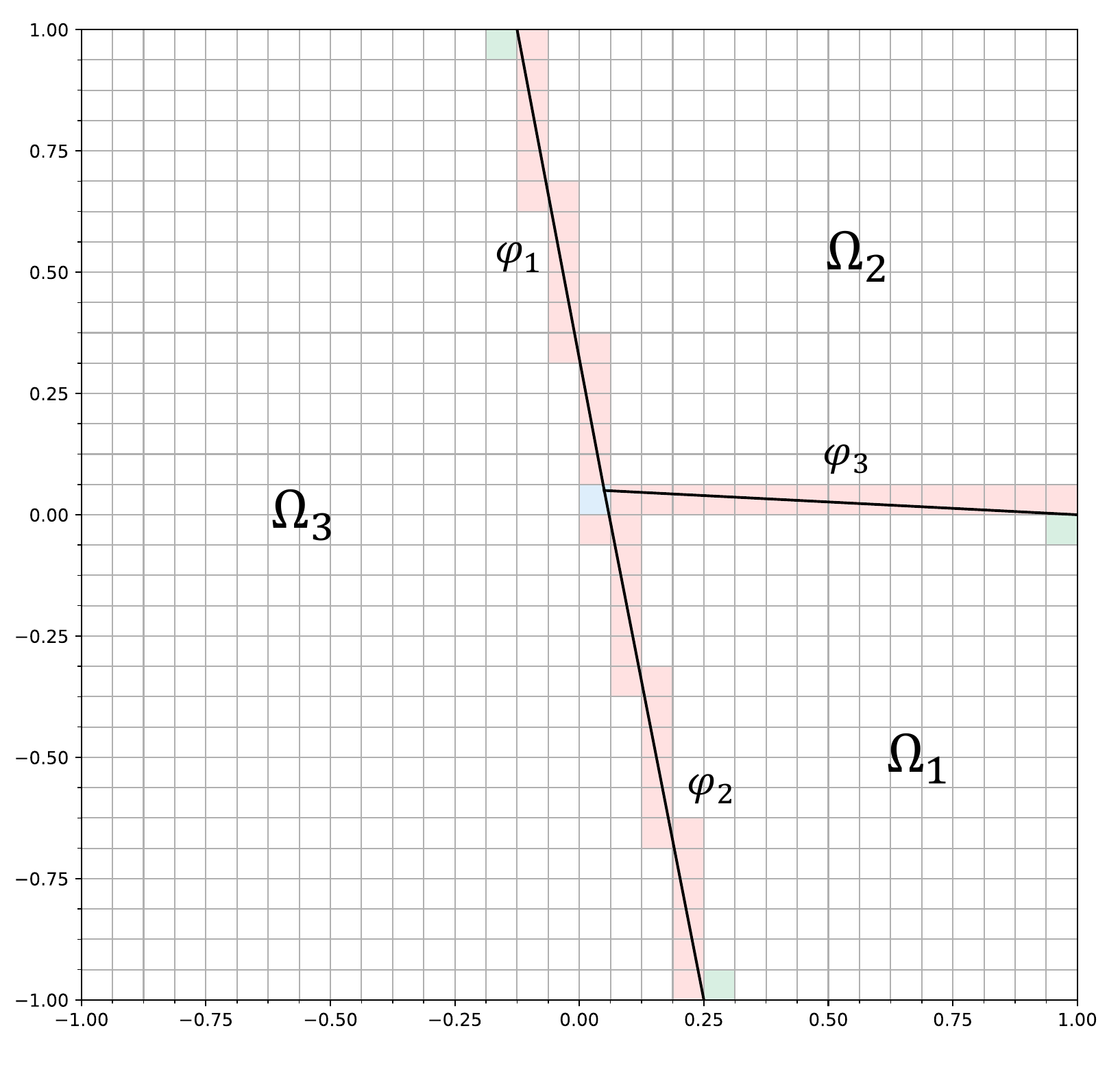}
  \includegraphics[width=.49\textwidth]{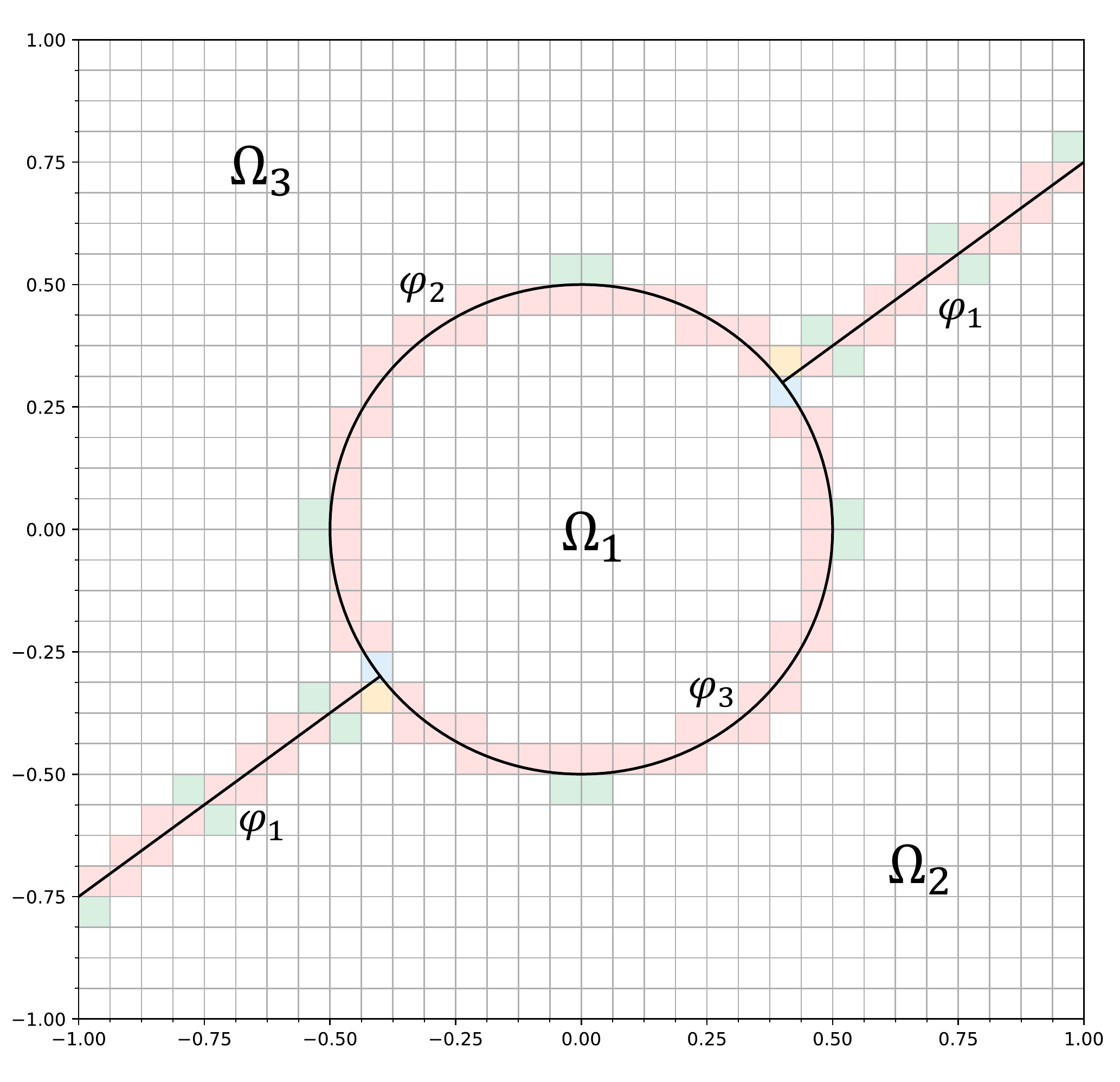}
  \caption{Geometry of the interface and element classification for Example 1 (left) and Example 2 (right)}
  \label{fig:Global Element classification of Numerical Examples}
\end{figure}

In this example, we consider a square domain $\Omega = (-1,1)^2$ separated by three straight-line interfaces. The interfaces are defined through the following level set functions:
\begin{equation*}
\begin{split}
&\varphi_1(x,y)=\frac{38}{7}x+y-\frac{9}{28},\\ 
&\varphi_2(x,y)=5.25x+y-0.3125, \\
&\varphi_3(x,y)=\frac{1}{19}x+y-\frac{1}{19}.
\end{split}
\end{equation*}
The exact solution of this problem is set to be
\begin{equation*}
\begin{split}
&u_1(x,y)=\dfrac{1}{\beta_1}\sin[(\frac{1}{19}x+y-\frac{1}{19})(5.25x+y-0.3125)],\\ 
&u_2(x,y)=\dfrac{1}{\beta_2}\sin[(\frac{1}{19}x+y-\frac{1}{19})(\frac{38}{7}x+y-\frac{9}{28})], \\
&u_3(x,y)=\dfrac{1}{\beta_3}\sin[(\frac{38}{7}x+y-\frac{9}{28})(5.25x+y-0.3125)/10].
\end{split}
\end{equation*}
The geometry of the interface and the element classifications are displayed in the left plot of Figure \ref{fig:Global Element classification of Numerical Examples}.  

We first consider the case with a moderate jump $(\beta_1,\beta_2,\beta_3) = (10,1,100)$. Interpolation errors in $L^2$, and semi-$H^1$ norms are reported in Table \ref{tab:1}. The errors and convergence rates of the PPIFEM solutions are reported in Table \ref{tab:SPP_Ex1_10_1_100}. It can be observed that both the interpolations and PPIFEM solutions converge with optimal rates in all three norms. For comparison, we also include the numerical results of Galerkin bilinear IFEM solutions without partial penalty in Table \ref{tab:IFEM_Ex1_10_1_100}. Convergence in the $L^\infty$ norm is not optimal. 

\begin{table}[H]
\centering
{
\begin{tabular}{ccccc}
	\toprule
	$N$ &$\|\mathcal{I}_hu-u\|_{L^2}$ & order&$|\mathcal{I}_hu-u|_{H^1}$&order\\
	\midrule
	$16$&$3.05\times 10^{-2}$&&$6.05\times 10^{-1}$&\\
	$32$&$7.71\times 10^{-3}$&1.98&$3.02\times 10^{-1}$&1.00\\
	$64$&$1.93\times 10^{-3}$&2.00&$1.51\times 10^{-1}$&1.00\\
	$128$&$4.84\times 10^{-4}$&2.00&$7.52\times 10^{-2}$&1.00\\
	$256$&$1.21\times 10^{-4}$&2.00&$3.76\times 10^{-2}$&1.00\\
	$512$&$3.02\times 10^{-5}$&2.00&$1.88\times 10^{-2}$&1.00\\
	\bottomrule
\end{tabular}
}
	\caption{Interpolation errors and convergences of Example 1 with coefficients (10,1,100)}
	\label{tab:1}
\end{table}

\begin{table}[H]
\resizebox{\textwidth}{!}{
\begin{tabular}{ccccccc}
	\toprule
	$N$ & $\|u_h-u\|_{L^\infty}$ & order&$\|u_h-u\|_{L^2}$ & order&$|u_h-u|_{H^1}$&order\\
	\midrule
	$16$&  $2.70\times 10^{-2}$&      &$2.81\times 10^{-2}$&       &$6.05\times 10^{-1}$&\\
	$32$&  $6.92\times 10^{-3}$&1.96&$7.10\times 10^{-3}$&1.98&$3.02\times 10^{-1}$&1.00\\
	$64$&  $1.75\times 10^{-3}$&1.99&$1.78\times 10^{-3}$&2.00&$1.51\times 10^{-1}$&1.00\\
	$128$&$4.38\times 10^{-4}$&2.00&$4.45\times 10^{-4}$&2.00&$7.52\times 10^{-2}$&1.00\\
	$256$&$1.09\times 10^{-4}$&2.00&$1.11\times 10^{-4}$&2.00&$3.76\times 10^{-2}$&1.00\\
	$512$&$2.74\times 10^{-5}$&2.00&$2.78\times 10^{-5}$&2.00&$1.88\times 10^{-2}$&1.00\\
	\bottomrule
\end{tabular}
}
\caption{PPIFEM errors in Example 1 with coefficients (10,1,100)}
	\label{tab:SPP_Ex1_10_1_100}
\end{table}

\begin{table}[H]
\resizebox{\textwidth}{!}{
\begin{tabular}{ccccccc}
	\toprule
	$N$ & $\|u_h-u\|_{L^\infty}$ & order&$\|u_h-u\|_{L^2}$ & order&$|u_h-u|_{H^1}$&order\\
	\midrule
	$16$&  $2.79\times 10^{-2}$&       &$2.83\times 10^{-2}$&       &$6.08\times 10^{-1}$&\\
	$32$&  $8.28\times 10^{-3}$&1.75&$7.16\times 10^{-3}$&1.98&$3.04\times 10^{-1}$&1.00\\
	$64$&  $2.44\times 10^{-3}$&1.76&$1.78\times 10^{-3}$&2.00&$1.51\times 10^{-1}$&1.01\\
	$128$&$7.61\times 10^{-4}$&1.68&$4.46\times 10^{-4}$&2.00&$7.55\times 10^{-2}$&1.00\\
	$256$&$3.04\times 10^{-4}$&1.33&$1.11\times 10^{-4}$&2.00&$3.77\times 10^{-2}$&1.00\\
	$512$&$1.32\times 10^{-5}$&1.20&$2.79\times 10^{-5}$&2.00&$1.89\times 10^{-2}$&1.00\\
	\bottomrule
\end{tabular}
}
\caption{IFEM errors in Example 1 with coefficients (10,1,100)}
	\label{tab:IFEM_Ex1_10_1_100}
\end{table}

Next we consider a larger coefficient jump, i.e., $(\beta_1,\beta_2,\beta_3) = (100,10000,1)$. The errors of PPIFEM and Galerkin IFEM are reported in Tables \ref{tab:SPP_Ex1_100_10000_1} and \ref{tab:IFEM_Ex1_100_10000_1}, respectively. It can be observed that the PPIFEM scheme is robust for this large contrast case in all three norms. Again, the convergence of $L^\infty$ norm for Galerkin IFEM is suboptimal. The magnitudes of the errors are significant larger than the errors of PPIFEM solutions.

\begin{table}[H]
\resizebox{\textwidth}{!}{
\begin{tabular}{ccccccc}
	\toprule
	$N$ & $\|u_h-u\|_{L^\infty}$ & order&$\|u_h-u\|_{L^2}$ & order&$|u_h-u|_{H^1}$&order\\
	\midrule
	$16$&  $8.87\times 10^{-3}$&      &$2.81\times 10^{-2}$&       &$6.44\times 10^{-1}$&\\
	$32$&  $2.93\times 10^{-3}$&1.60&$7.10\times 10^{-3}$&1.98&$3.25\times 10^{-1}$&0.99\\
	$64$&  $7.30\times 10^{-4}$&2.00&$1.78\times 10^{-3}$&2.00&$1.63\times 10^{-1}$&1.00\\
	$128$&$1.81\times 10^{-4}$&2.01&$4.45\times 10^{-4}$&2.00&$8.16\times 10^{-2}$&1.00\\
	$256$&$4.43\times 10^{-5}$&2.04&$1.11\times 10^{-4}$&2.00&$4.08\times 10^{-2}$&1.00\\
	$512$&$1.16\times 10^{-5}$&1.93&$2.78\times 10^{-5}$&2.00&$2.04\times 10^{-2}$&1.00\\
	\bottomrule
\end{tabular}
}
\caption{PPIFEM errors in Example 1 with coefficients (100,10000,1)}
	\label{tab:SPP_Ex1_100_10000_1}
\end{table}

\begin{table}[H]
\resizebox{\textwidth}{!}{
\begin{tabular}{ccccccc}
	\toprule
	$N$ & $\|u_h-u\|_{L^\infty}$ & order&$\|u_h-u\|_{L^2}$ & order&$|u_h-u|_{H^1}$&order\\
	\midrule
	$16$&  $9.27\times 10^{-3}$&      &$2.80\times 10^{-2}$&       &$6.43\times 10^{-1}$&\\
	$32$&  $2.84\times 10^{-3}$&1.71&$7.10\times 10^{-3}$&1.98&$3.25\times 10^{-1}$&0.98\\
	$64$&  $8.51\times 10^{-4}$&1.74&$1.78\times 10^{-3}$&2.00&$1.63\times 10^{-1}$&1.00\\
	$128$&$4.54\times 10^{-4}$&0.91&$4.45\times 10^{-4}$&2.00&$8.16\times 10^{-2}$&1.00\\
	$256$&$2.35\times 10^{-4}$&0.95&$1.11\times 10^{-4}$&2.00&$4.08\times 10^{-2}$&1.00\\
	$512$&$1.19\times 10^{-4}$&0.98&$2.79\times 10^{-5}$&2.00&$2.04\times 10^{-2}$&1.00\\
	\bottomrule
\end{tabular}
}
\caption{IFEM Errors in Example 1 with Coefficients (100,10000,1)}
	\label{tab:IFEM_Ex1_100_10000_1}
\end{table}


\subsection{Example 2}
\begin{figure}[htbp]
  \includegraphics[width=.49\textwidth]{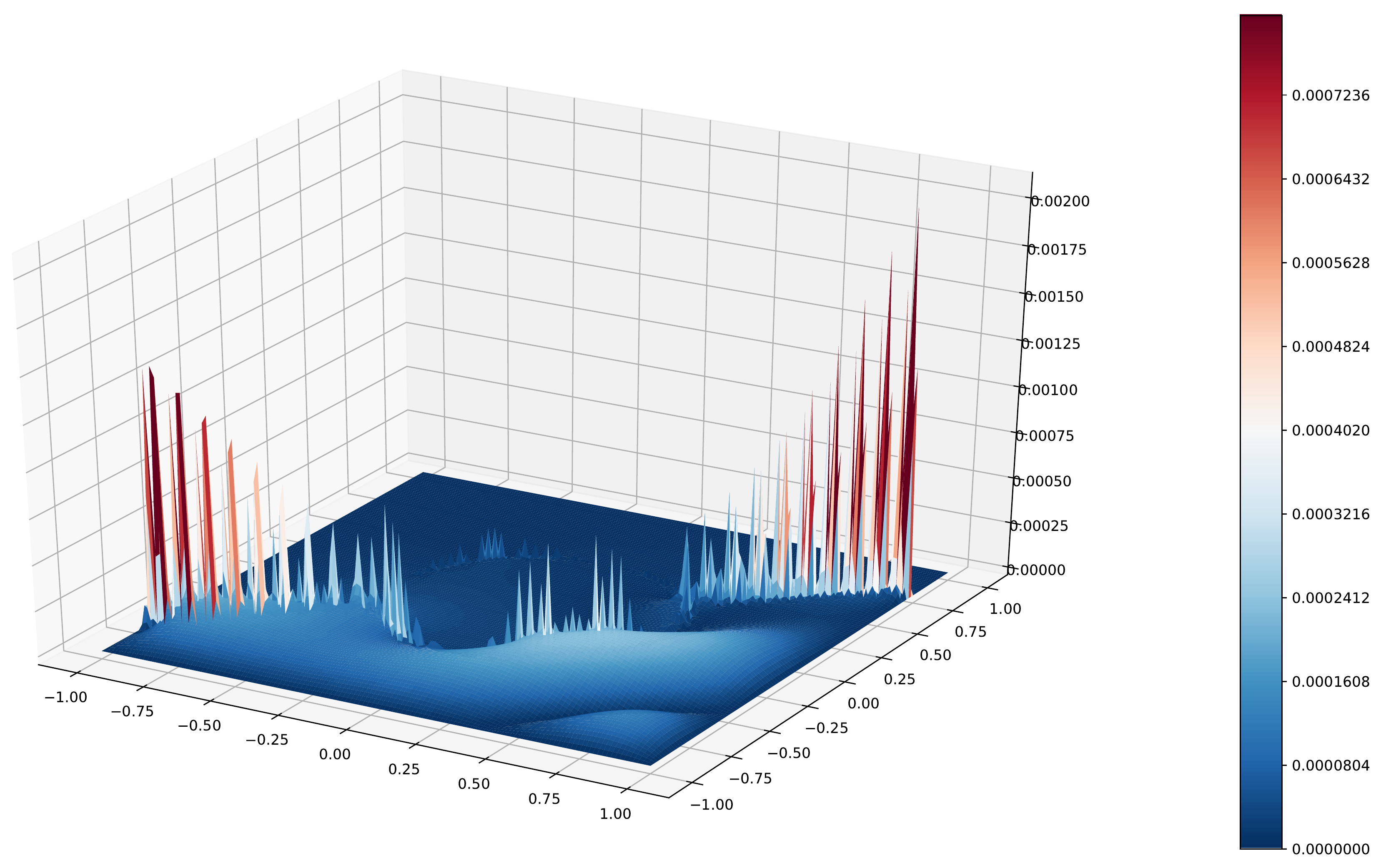}
  \includegraphics[width=.49\textwidth]{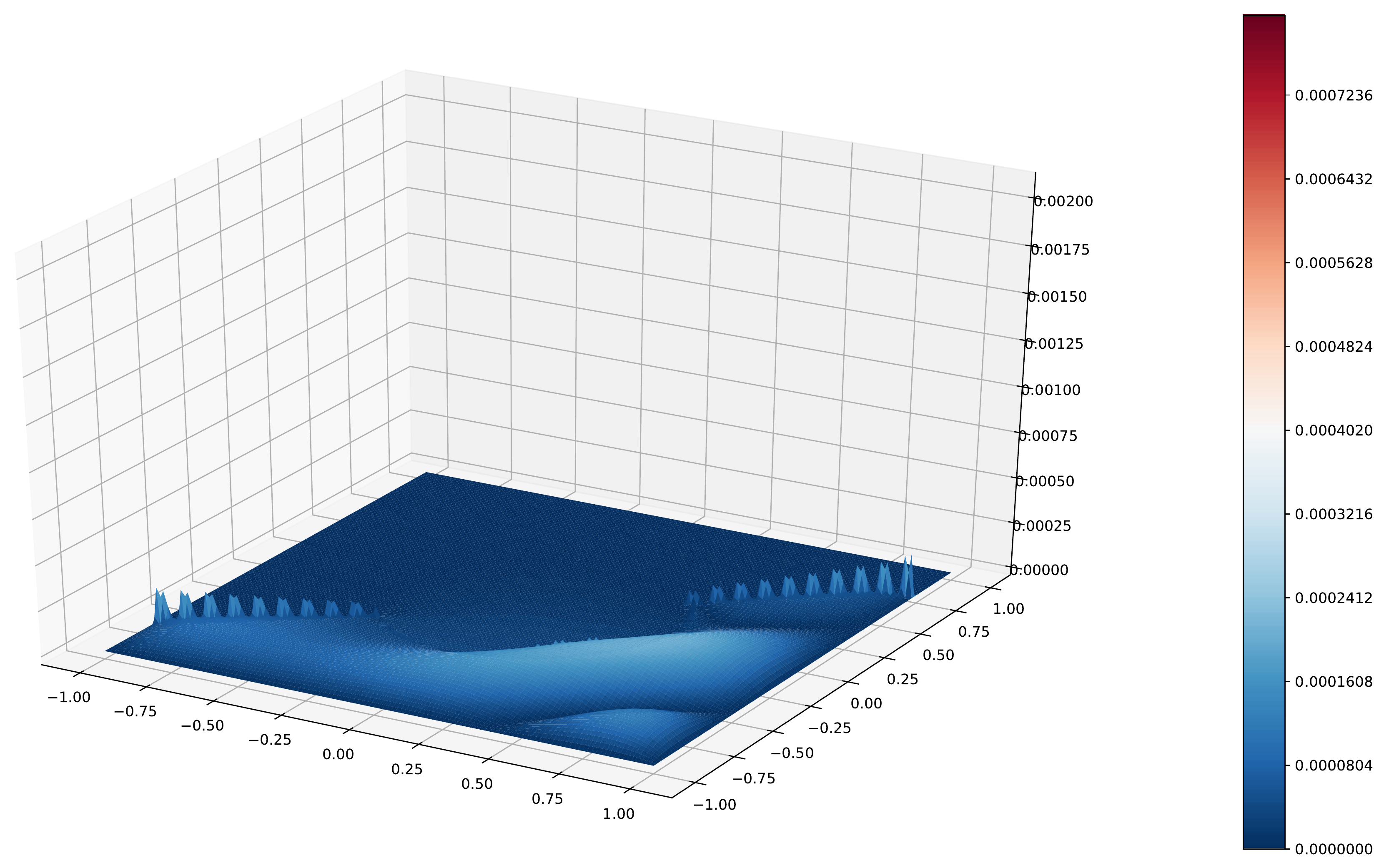}
  \caption{Comparison of error surfaces of classical IFEM and PPIFEM for Example 2.}
  \label{fig:difference}
\end{figure}

In this example, we test our numerical scheme for curved interfaces. In particular, the interfaces consist of a circle and a straight line, which are defined by the following level set functions 
\begin{equation*}
\begin{split}
&\varphi_1(x,y)=3x-4y, \\
&\varphi_2(x,y)=x^2+y^2-0.25,\\ 
&\varphi_3(x,y)=-x^2-y^2+0.25.
\end{split}
\end{equation*}
The exact solution is set to be
\begin{equation*}
\begin{split}
&u_1(x,y)=\dfrac{1}{\beta_1}((x^2+y^2)^{1.5}-0.125),\\ 
&u_2(x,y)=\dfrac{1}{\beta_2}(x^2+y^2-0.25)\sin(3x-4y), \\
&u_3(x,y)=\dfrac{1}{\beta_3}(3x-4y)\ln(x^2+y^2+0.75).
\end{split}
\end{equation*}
In the first case, we choose the coefficient to be $(\beta_1,\beta_2,\beta_3) = (10,1,100)$.  The error of the interpolation operator is reported in Table \ref{tab:2}. Errors of symmetric PPIFEM and Galerkin IFEM are reported in Tables \ref{tab:SPP_Ex2_10_1_100} and \ref{tab:IFEM_Ex2_10_1_100}, respectively. From these tables, we can see that the PPIFEM converge optimally in $L^\infty$, $L^2$ and $H^1$ norms. The Galerkin IFE solutions have suboptimal convergence in $L^\infty$ norm. Also, as the mesh size decreases, there seems to be an order deterioration in $L^2$ and $H^1$ norms for the Galerkin IFEM.  In Figure  \ref{fig:difference}, we compare the error surfaces of classical IFEM and the PPIFEM. This clearly demonstrates that PPIFEM outperforms the classical IFEM in terms of the accuracy around the interface.



\begin{table}[H]
\centering
{
\begin{tabular}{ccccc}
	\toprule
	$N$ &$\|\mathcal{I}_hu-u\|_{L^2}$ & order&$|\mathcal{I}_hu-u|_{H^1}$&order\\
	\midrule
	$16$&$2.46\times 10^{-2}$&&$5.27\times 10^{-1}$&\\
	$32$&$6.25\times 10^{-3}$&1.98&$2.63\times 10^{-1}$&1.00\\
	$64$&$1.57\times 10^{-3}$&1.99&$1.31\times 10^{-1}$&1.00\\
	$128$&$3.93\times 10^{-4}$&2.00&$6.57\times 10^{-2}$&1.00\\
	$256$&$9.84\times 10^{-5}$&2.00&$3.28\times 10^{-2}$&1.00\\
	$512$&$2.46\times 10^{-5}$&2.00&$1.64\times 10^{-2}$&1.00\\
	\bottomrule
\end{tabular}
}
	\caption{Interpolation errors and convergences of Example 2 with coefficients (10,1,100)}
	\label{tab:2}
\end{table}

\begin{table}[H]
\centering
\resizebox{\textwidth}{!}{
\begin{tabular}{ccccccc}
	\toprule
	$N$ & $\|u_h-u\|_{L^\infty}$ & order&$\|u_h-u\|_{L^2}$ & order&$|u_h-u|_{H^1}$&order\\
	\midrule
	$16$&$1.28\times 10^{-2}$&        &$2.20\times 10^{-2}$&&$5.26\times 10^{-1}$&\\
	$32$&$3.41\times 10^{-3}$&1.91&$5.55\times 10^{-3}$&1.99&$2.63\times 10^{-1}$&1.00\\
	$64$&$9.33\times 10^{-4}$&1.87&$1.39\times 10^{-3}$&2.00&$1.31\times 10^{-1}$&1.00\\
	$128$&$2.43\times 10^{-4}$&1.94&$3.48\times 10^{-4}$&2.00&$6.57\times 10^{-2}$&1.00\\
	$256$&$5.84\times 10^{-5}$&2.06&$8.71\times 10^{-5}$&2.00&$3.28\times 10^{-2}$&1.00\\
	$512$&$1.47\times 10^{-5}$&1.99&$2.18\times 10^{-5}$&2.00&$1.64\times 10^{-2}$&1.00\\
	\bottomrule
\end{tabular}
}
	\caption{PPIFEM Errors in Example 2 with Coefficients (10,1,100)}
	\label{tab:SPP_Ex2_10_1_100}
\end{table}

\begin{table}[H]
\centering
\resizebox{\textwidth}{!}{
\begin{tabular}{ccccccc}
	\toprule
	$N$ & $\|u_h-u\|_{L^\infty}$ & order&$\|u_h-u\|_{L^2}$ & order&$|u_h-u|_{H^1}$&order\\
	\midrule
	$16$&$4.75\times 10^{-2}$&        &$2.38\times 10^{-2}$&&$5.42\times 10^{-1}$&\\
	$32$&$1.49\times 10^{-2}$&1.67&$5.78\times 10^{-3}$&2.04&$2.68\times 10^{-1}$&1.01\\
	$64$&$4.82\times 10^{-3}$&1.63&$1.43\times 10^{-3}$&2.02&$1.34\times 10^{-1}$&1.01\\
	$128$&$2.09\times 10^{-3}$&1.21&$3.56\times 10^{-4}$&2.01&$6.70\times 10^{-2}$&1.00\\
	$256$&$9.64\times 10^{-4}$&1.11&$9.04\times 10^{-5}$&1.98&$3.39\times 10^{-2}$&0.98\\
	$512$&$4.74\times 10^{-4}$&1.02&$2.47\times 10^{-5}$&1.87&$1.75\times 10^{-2}$&0.96\\
	\bottomrule
\end{tabular}
}
	\caption{IFEM errors in Example 2 with coefficients (10,1,100)}
	\label{tab:IFEM_Ex2_10_1_100}
\end{table}

We also consider a large coefficient contrast case, $(\beta_1,\beta_2,\beta_3) = (100000,100,10)$. The errors of symmetric PPIFEM and Galerkin IFEM are reported in Tables \ref{tab:SPP_Ex2_1000000_100_10} and \ref{tab:IFEM_Ex2_1000000_100_10}, respectively. The convergences are similar to above.
\begin{table}[H]
\centering
\resizebox{\textwidth}{!}{
\begin{tabular}{ccccccc}
	\toprule
	$N$ & $\|u_h-u\|_{L^\infty}$ & order&$\|u_h-u\|_{L^2}$ & order&$|u_h-u|_{H^1}$&order\\
	\midrule
	$16$&$1.42\times 10^{-2}$&        &$3.29\times 10^{-3}$&&$5.77\times 10^{-2}$&\\
	$32$&$3.42\times 10^{-3}$&2.05&$7.22\times 10^{-4}$&2.19&$2.79\times 10^{-2}$&1.05\\
	$64$&$1.32\times 10^{-3}$&1.38&$2.11\times 10^{-4}$&1.78&$1.34\times 10^{-2}$&1.06\\
	$128$&$3.06\times 10^{-4}$&2.11&$3.60\times 10^{-5}$&2.55&$6.03\times 10^{-3}$&1.15\\
	$256$&$1.05\times 10^{-4}$&1.54&$8.60\times 10^{-6}$&2.07&$2.92\times 10^{-3}$&1.04\\
	$512$&$2.66\times 10^{-5}$&1.98&$2.07\times 10^{-6}$&2.06&$1.43\times 10^{-3}$&1.03\\
	\bottomrule
\end{tabular}
}
	\caption{PPIFEM Errors in Example 2 with coefficients (1000000,100,10)}
	\label{tab:SPP_Ex2_1000000_100_10}
\end{table}

\begin{table}[H]
\centering
\resizebox{\textwidth}{!}{
\begin{tabular}{ccccccc}
	\toprule
	$N$ & $\|u_h-u\|_{L^\infty}$ & order&$\|u_h-u\|_{L^2}$ & order&$|u_h-u|_{H^1}$&order\\
	\midrule
	$16$&$2.84\times 10^{-3}$&        &$2.07\times 10^{-3}$&&$4.59\times 10^{-2}$&\\
	$32$&$2.37\times 10^{-3}$&0.26&$5.66\times 10^{-4}$&1.87&$2.40\times 10^{-2}$&0.93\\
	$64$&$1.17\times 10^{-3}$&1.02&$1.85\times 10^{-4}$&1.61&$1.28\times 10^{-2}$&0.91\\
	$128$&$4.60\times 10^{-4}$&1.34&$3.77\times 10^{-5}$&2.30&$6.23\times 10^{-3}$&1.03\\
	$256$&$2.40\times 10^{-4}$&0.94&$1.02\times 10^{-5}$&1.88&$3.28\times 10^{-3}$&0.93\\
	$512$&$1.23\times 10^{-4}$&0.97&$3.15\times 10^{-6}$&1.70&$1.82\times 10^{-3}$&0.85\\
	\bottomrule
\end{tabular}
}
	\caption{IFEM errors in Example 2 with coefficients (1000000,100,10)}
	\label{tab:IFEM_Ex2_1000000_100_10}
\end{table}

\section{Conclusion}
We have developed a partially penalized immersed finite element method using bilinear polynomials for solving elliptic interface problems on multi-domain with triple junction points. The immersed finite element basis functions are constructed on interface elements with for all types of geometrical configurations. The local IFE space are enriched for handle nonhomogeneous flux jump conditions. Numerical results are provided to demonstrate the effectiveness of our method. In the future, we plan to carry on the a priori error estimates for the method and will also plan to extend the method to other interface problems.

\end{document}